\documentclass{svjour3}

\usepackage{amscd,amsfonts,amssymb,amsmath,latexsym,array,hhline,xcolor,graphicx}

\usepackage{latexsym}
\usepackage{times}

\newcommand\F{\mbox{I\kern-2pt F}}

\newcommand\cE{{\cal E}}

\newcommand\cF{{\cal F}}

\newcommand\e{{\varepsilon}}

\def\E{{\bf E}}
\def\P{{\bf P}}
\def\Chi{{\bf 1}}

\def\R{{\mathbb R}}
\def\bbr{{\mathbb R}}

\newcommand\fdem{$\Box$}
\newcommand\beq{\begin{equation}}
\newcommand\eeq{\end{equation}}
\newcommand\bea{\begin{eqnarray}}
\newcommand\eea{\end{eqnarray}}
\newcommand\bean{\begin{eqnarray*}}
\newcommand\eean{\end{eqnarray*}}

\begin{document}
\title{On  ruin probabilities with risky  investments  in a stock with stochastic volatility }

\author{ Anastasiya Ellanskaya 
\and 
 Yuri Kabanov}
\institute{  \at Lomonosov Moscow State University \and 
 \email{ellanskaya@gmail.com} \\
 \email{ykabanov@univ-fcomte.fr}    
 }
%
\date{Received: date / Accepted: date}

\titlerunning{On  ruin probabilities with risky  investments  }

\maketitle

\begin{abstract}
We investigate  the asymptotic of ruin probabilities when the company combines the  life- and non-life insurance businesses and invests its reserve into a risky asset with  stochastic volatility and drift driven by a two-state Markov process.  Using the technique of the implicit renewal theory 
we obtain the rate of convergence to zero of the ruin probabilities. 
\end{abstract}



 
   \section{Introduction}  
In the classical collective risk theory initiated by Filip Lundberg in 1903 and further developed by Harald Cram\'er in the  thirties it was usually assumed that an insurance company does not invest its reserve. This assumption leads to a good news: if claims 
are not heavy-tailed, then  the ruin probabilities are exponentially decreasing when the 
initial capital tends to infinity. 
For quite a long period of time studies of ruin probabilities did not suppose investments in spite that such an assumption  is clearly non realistic: in 
the modern world insurance companies operate in a financial environment and invests  their reserves into risky assets. Studies of mathematical models with risky investments started in the nineties. To the date there is an ample literature dedicated to the asymptotic  of  ruin probabilities  under various assumptions on the price processes  on the business activity.  
The fundamental discovery is  that the  risky investments change radically the asymptotic behavior of the ruin probabilities. For example, 
 in the papers Frolova, Kabanov and Pergamenshchikov  \cite{FrKP}, Kabanov and Pergamenshchikov \cite{KP}, and Kabanov and Pukhlyakov \cite{KPukh}  there were studied Lundberg type models with investments having, respectively,  the downward jumps (non-life insurance),  downward jumps (annuity model), and two-side jumps. 
  In the cited paper it was supposed that an insurance company  invests in a asset whose price process $S=(S_t)$ is a geometric Brownian motion (gBm) given by the equality $dS_t/S_t=a dt + \sigma W_t$. 
It was shown that  if the parameter $\beta:=2a/\sigma^2-1$ is strictly positive  then the ruin probability 
has the asymptotic behavior $C u^{-\beta}$ as the initial capital  $u$ tends to infinity. For the large volatility, when  $\beta\le 0$,  the ruin is imminent. 

A much more general setting, that of the generalized Ornstein--Uhlenbeck process, covering the case of investments with the  price following a geometric L\'evy process, was introduced by Paulsen in 1993, \cite{Paul-93},   and studied, using the techniques of implicit renewal theory, in the series of his papers, \cite{Paul-98,Paul-02,PG-97}. For the recent results 
in this direction and further references see \cite{KP2020}. 

To our knowledge,  the ruin problem with  risky investments in the asset with price with   
 stochastic volatility  was never studied. 
 In the present note we obtain the asymptotic of  ruin probabilities for the Lundberg--Cram\'er type model with investments and two-side jumps. We assume  that the price process 
 is a conditional geometric Brownian motion with  volatility and drift  modulated by a telegraph signal, that is, by a 
 two-state  Markov process, see the paper by Di Masi, Kabanov, and Runggaldier \cite{DKR}.  
 Such models are usually  referred to as the hidden Markov models or models with  regime switching.   
 In economics two-state  Markov processes  are used to model business cycles and we  believe that it is relevant to the situation with long  term investments typical in insurance. 

As in the paper   \cite{KP2020}, in our study we use the approach based on the recent progress in the implicit 
renewal theory. 

\section{The model}
We are given by a stochastic basis $(\Omega,\cF,{\bf F}=(\cF_t)_{t\ge 0},\P)$ with a Wiener process $W=(W_t)$, 
a  Poisson random measure $\pi(dt,dx)$ with
the compensator $\tilde \pi(dt,dx)=\Pi(dx)dt$, and a piecewise constant right-continuous  Markov process $\theta=(\theta_t)$ taking values in $\{0,1\}$ with 
transition intensity matrix $\Lambda=(\lambda^{ij})$ where $\lambda^{10}>0$, $\lambda^{01}>0$, $\lambda^{00}=-\lambda^{01}$, 
$\lambda^{11}=-\lambda^{10}$ and the initial value $\theta_0=i$. The $\sigma$-algebras generated by $W$, $\pi$, and $\theta$ are independent. 

 Let  $T_n$ be the successive jumps  
of the Poisson process $N$ with $N_t=p([0,t],\R)$ and let  $\tau_n$ be  the successive jumps  of $\theta$ with the convention that $T_0=0$ and $\tau_0=0$. 

Recall that the length between consecutive jumps of $\theta$ are independent random variables having exponential laws with parameters $\lambda^{01}$ and  $\lambda^{10}$ which are intensity of jumping from $0$ to $1$ and vice versa. So, 
if $\theta_0=0$, then the random variable $\tau_1$ has the exponential law with parameter 
$\lambda^{01}$, the random variable 
$\tau_2-\tau_1$ is independent on $\tau_1$ and has the exponential law with  parameter $\lambda^{10}$.

We consider the dynamics of the reserve $X=X^{u,i}$ of an insurance company where the business activity is as in the classical 
Lundberg--Cram\'er model but the reserve is fully invested  in a risky asset whose price 
$S$ follows a geometric Brownian motion with stochastic volatility and drift. We have that  
$$
dS_t=S_t (a_{\theta_t} dt+\sigma_{\theta_t}dW_t), \qquad S_0=1,   
$$ 
where $a_k\in \R$, $\sigma_k>0$, $k=0,1$.   In this case,
$X$ is of the form
\beq
 \label{risk}
 X_t=u+  \int_0^t  X_s  dR_s+ dP_t   
 \eeq
 where $R_t=R_0+a_{\theta_t}t+\sigma_{\theta_t} W_t$ 
 is the relative price  process  
 and 
\beq
 \label{r}
P_t=ct  + \int_0^t\int x  \pi(dt,dx)=ct+x*\pi_t
\eeq
is the process describing the business activity of the company. 

So, the  
model is described by the two-dimensional process $(X^{u,i},\theta^{i})$ where 
  $u>0$ is the initial capital and $i\in \{0,1\}$ is the initial value of $\theta$. 

We assume that the process $P$ is not increasing: otherwise the probability of ruin is zero. 

We also assume that ${\tilde \pi}(dt,dx)=dt\Pi(dx)$ with $\Pi(dx)=\alpha_1 F_1(dx)+\alpha_2F_2(dx)$ where 
$F_1(dx)$ is a probability distribution on $]-\infty,0[$ and $F_2(dx)$ is a probability distribution on $]0,\infty[$. In this
case the integral with respect to the jump measure is simply a difference of two independent  
compound Poisson processes with intensities $\alpha_1$, $\alpha_2$ for the jumps  downwards and upwards and whose absolute values have  the distributions $F_1(dx)$ and $F_2(dx)$ respectively. 

The solution of the linear equation (\ref{risk}) can be represented by the following stochastic version of the Cauchy formula  
\beq
\label{uY}
X_t^{u,i}:=\cE_t(R)(u-Y^i_t) 
\eeq 
where  
\beq
\label{Y_t}
Y_t^i:=-\int_{[0,t]} \cE^{-1}_{s}(R)dP_s=-\int_{[0,t]}  e^{-V_{s}}dP_s.  
\eeq
The  log price process $V=\ln \cE(R)$, that is  $V_s=\sigma_{\theta_s}W_s+(a_{\theta_s}-(1/2)\sigma_{\theta_s}^2)s$.  
 The process $Y^i$ does not depend on $u$. It will play the central role in our analysis.

 Let $\tau^{u,i}:=\inf \{t:\ X^{u,i}_t\le 0\}$ (the instant of ruin),
 $\Psi_i(u):=P(\tau^{u,i}<\infty)$ (the ruin probability),
 and $\Phi_i(u):=1-\Psi_i(u)$
 (the survival probability). 
It is easily seen that   $ \tau^{u,i}= \inf \{t\ge 0:\ Y^i_{t} \ge u\}$.

Recall that the constant parameter values $a=0$, $\sigma=0$, correspond to the 
Lundberg--Cram\'er model for which  the process $X^{u,i}=u+P_t$.  In the actuarial literature the compound Poisson process $P$   is usually written  in the form 
\beq
\label{r0}
P_t= ct  - \sum_{k=1}^{N_t}\xi_k
\eeq
where  either $\xi_k\ge 0$, $c> 0$ (i.e. $F_2=0$ --- jumps only downwards --- the case of non-life insurance)  or $\xi_k\le 0$, $c< 0$ (i.e. $F_1=0$ --- jumps only upwards --- the case of life insurance or annuity payments).  Models with both kinds of jumps 
are frequently  considered in the modern literature, see, e.g., \cite{AGY} and references therein.  
 For the classical models  with a positive average trend and $F$
having a ``non-heavy" tail, the Lundberg inequality asserts that the ruin probability decreases
exponentially as the initial capital of the company $u$ tends to infinity.
For the exponentially distributed claims the ruin
probability admits an explicit expression, see \cite{Asm},  Ch. IV.3b, or
\cite{Gr}, Section 1.1. 

For the models with exponentially distributed jumps with the investment in a risky asset with price following a geometric Brownian motion with the drift coefficient 
$a$ and the volatility $\sigma>0$ the answer is completely 
different: if  $2a/\sigma^2-1>0$ the ruin probability as a function of the initial capital $u$ behaves as $C u^{1-2a/\sigma^2}$ as $u\to \infty$. In the case where  
$2a/\sigma^2-1\le 0$ the ruin happens with probability one, see \cite{FrKP}, \cite{KP}, \cite{PZ}, \cite{KPukh}.  

To formulate our result for the model where the volatility and drift are  modulated by a telegraph process we introduce some notations. 

Define the  random variable $M_1:=e^{-V_{\tau_2}}$ and consider on $\R_+$ the moment generating function $f:q\mapsto \E M_1^q$.   
Put $\beta_i:=2a_i/\sigma_i^2-1$, $i=0,1$. 
It is easily seen that 
$$
f(q):=\E M_1^q= \E e^{-qV_{\tau_{2}}} =\E e^{-qV_{\tau_{1}}}  \E e^{-q(V_{\tau_{2}}-V_{\tau_1})}=f_0(q)f_1(q)
$$
where 
\bea
\label{0}
f_0(q)&:=&\frac{\lambda^{01}}{\lambda^{01}+q(a_0-\sigma_0^2/2)-q^2\sigma_0^2/2}=
\frac{\lambda^{01}}{\lambda^{01}+(1/2)\sigma_0^2q(\beta_0-q)},  \\ 
\label{1} 
f_1(q)&:=&\frac{\lambda^{10}}{\lambda^{10}+q(a-\sigma_1^2/2)-q^2\sigma_1^2/2}=
\frac{\lambda^{10}}{\lambda^{10}+(1/2)\sigma_1^2q(\beta_1-q)}
\eea
on the set ${\rm dom}\, f:=\{q\in \R_+\colon \E M_1^q<\infty\}$.  

Note that $f''(q)=q(q+1)\E V_{\tau_2}^2 e^{-qV_{\tau_{2}}}>0$, i.e. $f$ is strictly convex on its effective domain.  

 Suppose that $0<\beta_0< \beta_1$. Then $f_0(\beta_0)=1$, $f_1(\beta_0)<1$, hence, $f(\beta_0)<1$.  If $\beta_1\in {\rm dom}\, f$, then $f_0(\beta_1)>1$, $f_1(\beta_1)=1$,  and $f(\beta_1)>1$. If $\beta_1\notin {\rm dom}\, f$, then there is $\bar \beta\in ]\beta_0,\beta_1]$ such that $f_0(q)\to \infty$ as $q\uparrow \bar \beta$. Since $f_0$ is bounded away from zero on the interval $[\beta_0,\beta_1]$ we have that also
 $f(q)\to \infty$ as $q\uparrow \bar \beta$.  It follows that in all cases the convex function $f$ with $f(0)=1$ crosses the level $1$ also on the interval $ ]\beta_0,\beta_1[$, so  that there is a unique $\beta\in ]\beta_0,\beta_1[$ such that $f(\beta)=1$. In other words  $\beta$
 is the root of  the equation $f(q)=1$ which can be written in the detailed form as 
\beq
\label{root}
\sigma_0^2\sigma_1^2q(\beta_0-q)(\beta_1-q)+2\sigma_0^2(\beta_0-q)\lambda^{10}+2\sigma_1^2(\beta_1-q)\lambda^{01}=0.
\eeq   

It follows that 
\bea
\label{01}
(1/2)\sigma_0^2\beta(\beta_0-\beta)+\lambda^{01}=-\frac{\sigma_0^2(\beta_0-\beta)}{\sigma_1^2(\beta_1-\beta)}\lambda ^{10}>0,\\
\label{10} 
(1/2)\sigma_1^2\beta(\beta_1-\beta)+\lambda^{10}=-\frac{\sigma_1^2(\beta_1-\beta)}{\sigma_0^2(\beta_0-\beta)}\lambda^{01}>0.
\eea

\smallskip

Our main result is the following theorem. 

\begin{theorem}  
\label{main}
Suppose that $0<\beta_0< \beta_1$.  Let  $\beta\in ]\beta_0,\beta_1[$ be the solution of the equation   (\ref{root}). 
Suppose that $\Pi(|x|^\beta):=\int |x|^\beta \Pi(dx) <\infty$. 
Then for $i=0,1$
$$
0<\liminf_{u\to \infty} u^\beta \Psi_i(u)\le \limsup_{u\to \infty}u^\beta \Psi_i(u)<\infty. 
$$

  \end{theorem}
  
The proof  uses the techniques of the implicit renewal theory.  
 To apply  it, we verify that the random variables $
Q^i=Q^i_1:=- e^{-V^i}\cdot P_{\tau_2} $, $i=0,1$,   
belong to $L^\beta(\Omega)$, the process  $Y^i$ has at infinity a finite limit $Y^i_\infty$ which is a random variable unbounded from above and  has  the same law as $Q^i+MY^ i_{\infty}$ where 
$M=M_1$.  We already know that $\E M^\beta=1$ and $\E M^\beta+\e>0$ for some 
$\e>0$, and, of course, the law of $\ln M$ is not arithmetic. Moreover, we prove that 
$\bar G_i(u)\le \Psi_i (u)\le  C \bar G_i(u)$ where a constant $C>0$ and $\bar G_i(u)=\P(Y^i>u)$, Lemma \ref{G-Paulsen}. 

With such facts Theorem \ref{main}  follows from Theorem \ref{Le.3.2}  below which is the Kesten--Goldie theorem, see Th. 4.1 in \cite{Go91}, augmented by a statement on strict positivity of $C_+$ due to Guivarc'h and Le Page, \cite {BD} (for a simpler proof of the latter  see the paper \cite{BD}  by  Buraczewski and Damek and an extended discussion in Kabanov and Pergamenshchikov, \cite{KP2020}). 
\begin{theorem} 
\label{Le.3.2} 
Let $Y_{\infty}$ has the same law as $Q+MY_{\infty}$ where $M>0$.  
Suppose that $(M,Q)$ is such that the law of   $\ln\,M$ 
is non-arithmetic and, for  some $\beta>0$,
\begin{align}\label{3.3}
\E M^\beta=1, \ \ \ \E M^\beta\,(\ln\,M)^+<\infty, \  \  \ 
\E |Q|^\beta <\infty. 
\end{align}
Then 
\bean
&&\lim_{u\to\infty}\,u^\beta\,\P(Y_{\infty} >u)=C_+<\infty,\\
&&\lim_{u\to\infty}\,u^\beta\,\P(Y_{\infty} <-u)=C_-<\infty,
\eean
where $C_++C_->0$. 

If the  random variable  $Y_{\infty} $ is unbounded from above, then $C_+>0$. 
\end{theorem}

\section{Study of the process $Y$}
We consider the case where 
$\theta_0=0$ and omit $i=0$ in the notation.  The case $\theta_0=1$ is treated similarly, with the same results.  

The following identity is obvious: 
$$
Y_{\tau_{2n}}=-\sum_{k=1}^n \prod_{j=1}^{k-1}e^{-(V_{\tau_{2j}}-V_{\tau_{2j-2}})}\int_{]\tau_{2k-2},\tau_{2k}]}e^{-(V_s-V_{\tau_{2k-2}})}dP_s.
$$
Using the abbreviations
$$
Q_k:=- \int_{]\tau_{2k-2},\tau_{2k}]}e^{-(V_s-V_{\tau_{2k-2}})}dP_s, \qquad M_j:=e^{-(V_{\tau_{2j}}-V_{\tau_{2j-2}})} 
$$
we rewrite it in a more transparent form as 
\beq
\label{QM}
Y_{\tau_{2n}}=Q_1+M_1Q_2+M_1M_2Q_3+...+M_1...M_{n-1}Q_n. 
\eeq

The random variables $(Q_k,M_k)$ have the same law and are  
 is independent on the $\sigma$-algebra $\sigma\{(M_1,Q_1),..., (M_{k-1},Q_{k-1})\}$. 

\begin{lemma} 
\label{L-int}
 Let $\beta\in ]\beta_0,\beta_1[$ be the root of the equation (\ref{root}). Then 
\beq
\label{Leb}
\E\int_0^{\tau_2} e^{-\beta V_s}ds<\infty, \qquad 
\E  \Big (\int_{0}^{\tau_2}e^{-V_s}ds\Big)^\beta<\infty.
\eeq
\end{lemma}
{\sl Proof.} Using the independence of $\tau_1$ and $V$ we get  
that 
$$
\E\int_0^{\tau_1} e^{-\beta V_s}ds =\E \int_0^{\tau_1} e^{(1/2)\sigma_0^2\beta(\beta-\beta_0)s}
ds=\lambda^{01} \int_0^\infty e^{-[(1/2)\sigma_0^2\beta(\beta_0-\beta)+\lambda^{01}]t}dt<\infty 
$$
in virtue of the inequality (\ref{01}).  
Also, 
$$
\E\int_{\tau_1}^{\tau_2} e^{-\beta V_s}ds =\E e^{- \beta V_{\tau_1}}\int_{\tau_1}^{\tau_2} e^{-\beta (V_s-V_{\tau_1})}ds= \E e^{-\beta V_{\tau_1}}\E \int_{\tau_1}^{\tau_2} e^{-\beta (V_s-V_{\tau_1})}ds. 
$$
Recall that $\E e^{-\beta V_{\tau_1}}=f_0(\beta)<\infty$. Taking into account that on the 
interval $[\tau_1,\tau_2]$ the process $V_s-V_{\tau_1}$ is a Wiener process with variance 
$\sigma_1^2$ and drift $(1/2)\sigma_1^2\beta_1$, we obtain that   
$$
\E \int_{\tau_1}^{\tau_2} e^{-\beta (V_s-V_{\tau_1})}dt=\lambda^{10}
 \int_0^\infty e^{-[(1/2)\sigma_1^2\beta(\beta_1-\beta)+\lambda^{10}]t}dt<\infty 
$$ 
since  both terms in the square brackets are strictly positive (and we have no need to use the representation in (\ref{10})).  This gives us the first property in (\ref{Leb}).

Take $p>1$ such  that $(\beta/p- \beta_0)>0$ and  define the constant $\kappa:=(1/2)\sigma_0^2\beta/p$ and martingale 
$L=(L_s)$ with $L_s=e^{(-\sigma_0\beta W_s-\kappa\beta s)/p}$. Then 
$$
 \Big (\int_0^{t}e^{-V_s}ds\Big)^\beta=\Big(\int_0^{t}e^{-\sigma_0 W_s-\kappa s}e^{\kappa s-(1/2)\sigma_0^2\beta_0 s}ds\Big)^\beta\le L^{*p}_t \Big(\int_0^{t}e^{\kappa s-(1/2)\sigma_0^2 \beta_0s}ds\Big)^\beta
$$
where $L^*_t:=\sup_{s\le t}L_s$.  By the Doob inequality 
$$
\E L^{*p}_t\le C_p \E L^p_t= C_p e^{(1/2)\sigma_0^2\beta^2(1-1/p)t}.
$$
where $C_p:=(p/(p-1))^p$. Integrating the exponential we get that 
$$
\Big(\int_0^{t}e^{\kappa s-(1/2)\sigma_0^2 \beta_0s}ds\Big)^\beta=
\Big(\int_0^{t}e^{(1/2)\sigma_0^2(\beta/p- \beta_0)s}ds\Big)^\beta\le (1/c_p) e^{c_pt}
$$
where  $c_p:=(1/2)\sigma_0^2\beta (\beta/p- \beta_0)$. 

Combining two estimates we get that 
$$
\E  \Big (\int_0^{t}e^{-V_s}ds\Big)^\beta\le (C_p/c_p) e^{(1/2)\sigma_0^2\beta (\beta- \beta_0)t}<\infty  
$$
and 
$$
\E  \Big (\int_0^{\tau_1}e^{-V_s}ds\Big)^\beta\le \lambda^{01} (C_p/c_p) \int_0^\infty e^{-[(1/2)\sigma_0^2\beta(\beta_0-\beta)+\lambda^{01}]t}dt<\infty 
$$
in virtue of (\ref{01}).  Using  again the factorization 
$$
 \E  \Big (\int_{\tau_1}^{\tau_2}e^{-V_s}ds\Big)^\beta=\E e^{-\beta V_{\tau_1}}\E \Big (\int_{\tau_1}^{\tau_2}e^{-(V_s-V_{\tau_1})}ds\Big)^\beta
$$
and observing that 
 $$
 \E  \Big (\int_{\tau_1}^{\tau_2}e^{-V_s}ds\Big)^\beta\le \lambda^{10}(C_p/c_p) \int_0^\infty e^{-[(1/2)\sigma_0^2\beta(\beta_1-\beta)+\lambda^{01}]t}dt<\infty.  
 $$
we get the second property in (\ref{Leb}). \fdem

\smallskip The following lemma that our conditions imply that  $Q_1\in L^\beta(\Omega)$. 

\begin{lemma}  
\label{Q1}
Suppose that $\Pi(|x|^{\beta}):=\int |x|^{\beta}\Pi(dx)<\infty$. Then  $\E |Q_1|^{\beta}<\infty$.
\end{lemma}
{\sl Proof.}
{\bf Case where  $\beta\le 1$}. The inequality  $(|x|+|y|)^\beta\le |x|^\beta+|y|^\beta$ allows us to consider the integrability of summands separately. 
In particular, for the jump component 
of the process $P$ whose jump measure $\pi(dt,dx)$ has the compensator $\tilde \pi(dt,dx)=\Pi(dx)dt$), we get  that 
$$
\E (e^{-V}x  *\pi_{\tau_2})^\beta\le \E e^{-\beta V}|x|^\beta  *\pi_{\tau_2}=\E e^{-\beta V}|x|^\beta  *\tilde \pi_{\tau_2}\le \Pi(|x|^\beta)\E\int_0^{\tau_2} e^{-\beta V_s}ds. 
$$
The claim follows from here and Lemma \ref{L-int}. 

{\bf Case where  $\beta> 1$.} Now we shall split summands using the elementary  inequality 
$$
(|x|+|y|)^\beta\le 2^{\beta-1}(|x|^\beta+|y|^\beta). 
$$
Because of Lemma \ref{L-int} we need to consider only the integral with respect to the 
jump component of $P$. 
Note that $e^{-V}|x|  *\tilde \pi_{\tau_2}<\infty$. Then 
$$
\E \big(e^{-V}|x| *\pi_{\tau_2}\big )^\beta\le  2^{\beta-1}\big (\E \big(e^{-V}|x|  *(\pi-\tilde \pi)_{\tau_2}\big)^\beta + 
\E\big( e^{-V}|x|  *\tilde \pi_{\tau_2}\big)^\beta\big).
$$
Due to (\ref{Leb})
$$
\E\big( e^{-V}|x|  *\tilde \pi_{\tau_2}\big)^\beta\le \big(\Pi(|x|) \big)^\beta \E  \Big (\int_{0}^{\tau_2}e^{-V_s}ds\Big)^\beta<\infty.
$$
Let $I_s:=e^{-V}|x|  *(\pi-\tilde \pi)_s$. According to the Novikov inequalities (with $\alpha=1$) 
the moment of the order $\beta>1$ of $I^*_t:=\sup_{s\le t}|I_s|$ admits the bound  
\bean
\E I^{*p}_{\tau_2} &\le &C_{\beta,1}\big(\E (e^{-V}|x|*\tilde \pi_{\tau_2} )^\beta + \E e^{-\beta V}|x|^\beta *\tilde \pi_{\tau_2}  \big)\\
&\le & C'_{\beta,1}\E \Big(\int_0^{\tau_2} e^{-V_s}ds\Big )^\beta+ \E
C'_{\beta,1}\int_0^{\tau_2} e^{-\beta V_s}ds\Big. 
\eean
where $C'_{\beta,1}:=C_{\beta,1}(\Pi(|x|))^\beta<\infty$,  $C''_{\beta,1}:=C_{\beta,1}\Pi(|x|^\beta)<\infty$ due to our assumption. 
The both integrals in the right-hand side  as we proved are finite.  \fdem

\begin{lemma} The process $Y$ has the following properties: 

$(i)$ $Y_t$ converges almost surely as $t\to \infty$ to a finite random variable 
$Y_\infty$.

$(ii)$ $Y_\infty=Q_1+M_1Y_{1,\infty}$ where $Y_{1,\infty}$ is a random variable independent on $(Q_1,M_1)$ and having the same law as $Y_{\infty}$. 

$(iii)$ $Y_{\infty}$ is unbounded from above. 
\end{lemma}
{\sl Proof.} $(i)$ Take $p\in ]0,\beta\wedge 1 [$. Then $r:=\ EM_1^p<1$ and, by virtue of Lemma \ref{Q1}, $\E |Q_1|^p<\infty$. It follows that 
$\E |Y_{\tau_{2n+2}}-Y_{\tau_{2n}}|^p=\E M_1^p\dots M_n^p Q_{n+1}^p=r^n \E |Q_1|^p$ and, therefore, 
$$
\E \Big(\sum_{n\ge 0} |Y_{\tau_{2n+2}}-Y_{\tau_{2n}}|\Big )^p\le \sum_{n\ge 0}  \E |Y_{\tau_{2n+2}}-Y_{\tau_{2n}}|^p<\infty.
$$ 
Thus, $\sum_n |Y_{\tau_{2n+2}}-Y_{\tau_{2n}}|<\infty$ a.s. implying that $Y_{\tau_{2n}}$ converges a.s. to some finite random variable we shall denote $Y_{\infty}$. 

Let $Y_t^*:=\sup_{s\le t}|Y_s|$.  Then  
$$
\E Y^{*p}_{\tau_2}
\le 
c^p  E\Big(\int_0^{\tau_2}e^{-V_t}dt\Big)^p +(\Pi(|x|))^p E\Big(\int_0^{\tau_2}e^{-pV_t}dt\Big)^p<\infty.
$$
Put 
$$
\Delta_n=\sup_{v\in [\tau_{2n},\tau_{2n+2}]}\Big\vert\int_{\tau_{2n}}^v e^{-V_s}dP_s\Big\vert.
$$
Then
$$
\E \Delta_n^p= \E M_1^p\dots M_n^p \sup_{v\in [\tau_{2n},\tau_{2n+2}]}\Big\vert\int_{\tau_{2n}}^v e^{-(V_{s}-V_{\tau_{2n}})}dP_s\Big\vert^p\le r^p\E Y^{*p}_{\tau_2}
$$
and, therefore, for any $\e>0$
$$
\sup_{n\ge 0}\P ( \Delta_n\ge \e)\le e^{-p}\sup_{n\ge 0}\E \Delta_n^p<\infty. 
$$
By the Borel--Cantelli lemma for all $\omega$ except a null-set $\Delta_n(\omega)\le \e$ for all $n\ge n(\omega)$. This implies that $Y_t$ converges a.s.  to the same limit as the sequence $Y_{\tau_n}$. 
\smallskip

$(ii)$  Rewriting (\ref{QM}) in the form 
$$
Y_{\tau_{2n}}=Q_1+M_1(Q_2+M_2Q_3+...+M_2...M_{n-1}Q_n)
$$
and observing that the sequence of random variables in the parentheses converges 
as to a random variable with same law as $Y_{\infty}$ and independent on $(Q_1,M_1)$ we get the needed assertion. 

$(iii)$ In virtue of $(ii)$ it is sufficient to check that  the set 
$\{Q_1\ge N,\ M_1\le 1/N\}$ is non-null whatever is $N\ge 0$. We consider several 
cases.
 
$1)$ $c<0$.   On the set $\{T_1>\tau_2\}$  we have 
$Q_1=-c\int_0^{\tau_2} e^{-V_s}ds$. For every $t>0$ the set 
$$
B_N(t):=\Big\{-c\int_0^t e^{-V_s}ds\ge N,\ e^{-V_t}\le 1/N\Big\}
$$ 
is non-null and so is $B_N(\tau_2)$.  
Note that  
$$
\Gamma_N:=\{Q_1\ge N,\ M_1\le 1/N,\ T_1>\tau_2\}=B_N(\tau_2)\cap \{T_1>\tau_2\}. 
$$ 
Thus, 
$$
\P (\Gamma_N):= \E I_{B_N(\tau_2)}I_{\{T_1>\tau_2\}}=\E I_{B_N(\tau_2)}e^{-(\alpha_1+\alpha_2)\tau_2}>0. 
$$

$2)$ $c\ge 0$. Since $P$ is not increasing, $\Pi(]-\infty,0[)>0$. Let  $r_N:=\delta + \ln N$. We consider the set  
$$
\Delta_N:= \{ -\delta+\gamma< V_s<\delta+\gamma s,\ \forall\,s\le r_N+1\}\cap \{ r_N\le \tau_2\le r_N+1 \}
$$ 
The set $\Delta_N$ is non-null and on this set $e^{-V_{\tau_2}}\le  e^{\delta -r_N} \le 1/N$ and 
$$
c\int_0^{\tau_2}e^{-V_s}ds\ge \int_0^{\tau_2}e^{-V_s}ds
$$  
Note that the set 
$$
\Big\{e^{-V}*x_{\{x< 0\}}\pi_{r_N}\ge C_N +N,\ x_{\{x >0\}}*\pi_{\tau_2}=0 \Big\}
$$
is non-null and so non-null is its intersection with $\Delta_N$. But this intersection  is a subset of the set $\{Q_1\ge N,\ M_1\le 1/N\}$. \fdem 

\begin{lemma}
\label{G-Paulsen} 
For all  $u>0$  
\beq
\label{Paulsen}
\bar G_i(u)\le\, \Psi_i(u)=\frac{\bar G_i(u)}{\E\big(\bar G_{\theta_{\tau^{u,i}}}(0) \vert \tau^{u,i}<\infty\big)}\le \frac{\bar G_i(u)}{\bar G_0(0)\wedge \bar G_1(0)},
\eeq  
where $\bar G_i(u):=\P(Y_\infty^i>u)$. 
\end{lemma}
\noindent
{\sl Proof.}
Let $\tau$ be an arbitrary stopping time with respect to the  filtration $(\cF^{P,R,\theta}_t)$. 
As the finite limit $Y^i_\infty$ exists,  the random variable 
$$
Y_{\tau,\infty}^i:=\begin{cases}
-\lim_{N\to \infty } 
\int_{]\tau,\tau+N]}\,e^{-(V_{s}-V_{\tau})}dP_{s},&\tau<\infty, \\
0, & \tau=\infty,
\end{cases} 
$$ 
is well defined.  On the set $\{\tau<\infty\}$
\beq
\label{YX}
Y_{\tau,\infty}^i=e^{V_\tau}(Y^i_{\infty}-Y^i_\tau)=X_{\tau}^u  +e^{V_\tau}(Y^i_\infty-u). 
\eeq   Let $\xi$ be a $\cF_{\tau}^{P,R,\theta}$-measurable random variable.  Note that  the conditional distribution of $Y_{\tau,\infty}^i$ given $(\tau,\xi,\theta_\tau)=(t,x,j)\in \bbr_+\times\bbr\times \{0,1\}$ is the same as the distribution of $Y_{\infty}^j$. It follows that   
$$
\P\big (
Y_{\tau,\infty}^i>\xi, \ 
\tau<\infty,\ \theta_\tau=j
\big) 
=\E \bar G_j(\xi)\, \Chi_{\{ \tau<\infty,\; \theta_\tau=j\}}.
$$
Thus, if $\P(\tau<\infty)>0$, then
$$
\P\big(Y_{\tau,\infty}^i >\xi, \ 
\tau<\infty
\big)
=\E\big(\bar G_{\theta_\tau}(\xi)\, \vert\, \tau<\infty\big)\,\P\big(\tau<\infty\big)\,.
$$
Noting that  $\Psi_i(u):=\P\big(\tau^{u,i}<\infty\big)\ge \P\big(Y^i_{\infty}>u\big)=\bar G_i(u)>0$,   we deduce from here using (\ref{YX}) that 
\begin{align*}
\bar G_i(u)&=
\P\big(
Y^i_{\infty}>u,\ \tau^{u,i}<\infty\big)=
\P\big(Y_{\tau^{u,i},\infty}^i>X_{\tau^{u,i}}^{u,i},\ 
\tau^{u,i}<\infty
\big)\\
&=\E\big(\bar G_{\theta_{\tau^{u,i}}}(X_{\tau^{u,i}}^{u,i})\, \vert\, \tau^{u,i}<\infty\big)\,\P\big(\tau^{u,i}<\infty\big)
\end{align*}
implying the equality in (\ref{Paulsen}). Also, 
$$
\E\big(\bar G_{\theta_{\tau^{u,i}}}(0) \vert \tau^{u,i}<\infty\big)=
\bar G_0(0)\P\big(\theta_{\tau^{u,i}}=0\vert \tau^{u,i}<\infty\big)+\bar G_1(0)\P\big(\theta_{\tau^{u,i}}=1\vert \tau^{u,i}<\infty\big). 
$$
The result follows since  
$X_{\tau^{u,i}}^{u,i}\le 0$ on   $\{\tau^{u,i}<\infty\}$ and,  in the case  with only upward jumps that is when 
$\Pi(]-\infty,0[)=0$, the process $X^{u,i}$ crosses zero in a continuous way, i.e.  $X_{\tau^{u,i}}^{u,i}= 0$ on this set. \fdem


\smallskip
{\bf Acknowledgement.} This work was supported by the Russian Science Foundation grant  
20-68-47030.  


 \end{document}